\documentclass[12pt]{amsart}
\usepackage{amssymb}
\usepackage{mathptmx}

\newtheorem{theorem}{Theorem}
\newtheorem{lemma}{Lemma}

\newtheorem{corollary}{Corollary}
\newtheorem{remark}{Remark}

\begin{document}

\title[]{On the Riesz Basis Property\ of the Eigen- and Associated Functions \ of
Periodic and Antiperiodic Sturm-Liouville Problems}

\author[]{A.A.~Shkalikov \ \, and
O.A.~Veliev}
\thanks{A.A.Shkalikov is supported by
Russian Foundation of Basic Research (project No 07-01-00283)
 and by INTAS (project No 05-1000008-7883).}

\date{}
\maketitle

\begin{abstract}
The paper deals with the Sturm-Liouville operator
\[
Ly=-y^{\prime\prime}+q(x)y,\qquad x\in\lbrack0,1],
\]
generated in the space $L_{2}=L_{2}[0,1]$ by periodic or antiperiodic boundary
conditions. Several theorems on Riesz basis property of the root functions of
the operator $L$ are proved. One of the main results is the following.
\textsl{Let $q$ belong to Sobolev space $W_{1}^{p}[0,1]$ with some integer
$p\geq0$ and satisfy the conditions $q^{(k)}(0)=q^{(k)}(1)=0$ for
$0\leq k\leq s-1$, where s}$\leq p.$ \textsl{Let the functions $Q$
and $S$ be defined
by the equalities $Q(x)=\int_{0}^{x}q(t)\,dt,\ S(x)=Q^{2}(x)$ and let $q_{n}%
,Q_{n},S_{n}$ be the Fourier coefficients of $q,Q,S$ with respect
to the trigonometric system $\{e^{2\pi inx}\}_{-\infty}^{\infty}$.
Assume that the sequence $q_{2n}-S_{2n}+2Q_{0}Q_{2n}$ decreases
not faster than the powers $n^{-s-2}$. Then the system of eigen
and associated functions of the operator $L$ generated by periodic
boundary conditions forms a Riesz basis in the space $L_{2}[0,1]$
(provided that the eigenfunctions are normalized) if and only if
the condition
\[
q_{2n}-S_{2n}+Q_{0}Q_{2n}\asymp
q_{-2n}-S_{-2n}+2Q_{0}Q_{-2n},\quad n>1,
\]
holds.}

\end{abstract}


\section{Introduction and preliminary results}

Denote by $L_{\alpha}$ the Sturm-Liouville operator
\[
L_{\alpha}=-y^{\prime\prime}+q(x)y,
\]
generated in the space $L_{2}[0,1]$ by quasi-periodic boundary conditions
\[
y(1)=e^{\pi i\alpha}y(0), \quad y'(1)=e^{\pi i\alpha}y'(0).
\]
Generally we assume that the potential $q(x)$ is a complex valued
Lebesque summable function and $\alpha\in\lbrack0,2)$. It is known
that the operator $L_{\alpha}$ is Birkhoff regular. Moreover, in
the case $\alpha\neq0,1$ it is strongly regular. In the latter
case (i.e. for $\alpha\neq0,1$) the normalized eigen- and
associated functions (or the root functions) of the operator
$L_{\alpha}$ form a Riesz basis (this result is proved
independently in \cite[ Ch.XIX]{DS}, \cite{Ke}, \cite{Mi} for
 strongly regular differential operators  of arbitrary order rather
than for the Sturm-Liouville operators). In the case when an
operator is regular but not strongly regular the root functions,
generally, do not form even usual basis. However, it is known
\cite{Sh1},
\cite{Sh2}, that they can be combined in pairs, so that the
corresponding 2-dimensional subspaces form a Riesz basis of
subspaces (see the definitions in \cite[Ch. 6]{GK}, for example).

The aim of this paper is to study in details the cases $\alpha=0$
and $\alpha=1$ which correspond to periodic and antiperiodic
problems, and to find necessary and sufficient conditions for the
Riesz basis property of the root functions. This problem is
important for the study of Sturm-Liouville operators with periodic
complex potential on the whole real line (the so-called Hill
operator) and it is of independent interest. To make acquaintance
with recent progress in the study of the Hill operator we refer
the readers to a paper by Djakov and Mitjagin
\cite{DM1} where other many references can be found.

 In 1996 at the seminar in MSU A.A.Shkalikov  formulated
the following result:
\textsl{Assume that $q(x)$ is a smooth potential, $q^{(k)}(0)=q^{(k)}(1)=0$
for $0\leq k\leq s-1$, and $q^{(s)}(0)\neq q^{(s)}(1)$. Then the
normalized root functions of the operator $L_{0}$ form a Riesz
basis in $L_{2}$.} The idea of the proof was the following: to
obtain the asymptotics of the eigenfunctions using the well-known
asymptotic of the solutions of the Sturm-Liouville equation with
smooth potentials and then  to check that the angle between the
eigenfunctions in 2-dimensional subspaces does not tend to zero.
Then the basis property follows from theorems of
\cite{Sh1}, \cite{Sh2}. Kerimov and Mamedov \cite{KM}
 obtained the rigorous
proof of this result in the case $q\in C^4[0,1], \ q(1)\ne q(0)$.
Actually,  this results remains valid for an arbitrary
 $s\ge 0$.   It will be obtained in Corollary~2 of
this paper.

Another more original approach is due to Dernek and Veliev
\cite{DV}. The result was obtained in terms of the Fourier
coefficients of the potential $q.$ Namely, they proved the
following result:
\textsl{Assume that for $\alpha =0$ or for $\alpha =1$  the conditions
\begin{gather}\label{1}
\lim_{n\rightarrow\infty}\frac{\ln\left\vert n\right\vert }{nq_{2n+\alpha}%
}=0,\\
\label{2}
q_{2n+\alpha}\asymp q_{-2n-\alpha},
\end{gather}
hold. Then the root functions of the operator $L_{\alpha}$,
generated by periodic or antiperiodic boundary conditions,
respectively, form a Riesz basis in $L_{2}[0,1]$. }

Makin \cite{M} essentially improved this result. Using another
method he proved that  the assertion on the Riesz basis property
remains valid  if condition \eqref{2} holds, but condition
\eqref{1} is replaced by a less restrictive one:
\begin{gather*}
q\in W_{1}^{s}[0,1],\quad q^{(k)}(0)=q^{(k)}(1),\quad\text{for}%
\ \,k=0,1,...,s-1,\\
\mid q_{2n+\alpha}\mid>c_{0}n^{-s-1},\quad\forall\ n\gg1\quad\text{with
some}\ \,c_{0}>0,
\end{gather*}
where $s$ is a nonnegative integer. Moreover, some
conditions which involve the absence of the Riesz
basis property were presented in \cite{M}. Some sharp
results on the absence  of the Riesz basis property
were obtain by Djakov and Mitjagin \cite[Theorem
71]{DM1}.

The results which we obtain in the present paper are
more general and cover all the previous ones, except
constructions in Theorem 71 of \cite{DM1} (see
Corollaries at the end of the paper).

Further we shall work only  with periodic problem. The changes
which have to be done for antiperiodic problem are obvious, and we
shall note on them at the end of the paper. Without loss of
generality, we always  assume that $q_{0}=\int_{0}
^{1}q(x)\,dx=0$. It is known \cite[Ch.1]{Na}, that the eigenvalues
of $L_{0}$ are grouped in pairs. It is convenient to distribute
them  in two sequences $\{\lambda_{n,1}\}_{1}^{\infty},\
\{\lambda _{n,2}\}_{1}^{\infty}$ such that each of these sequences
(see \cite{DV}, for example) obeys the asymptotic
\begin{equation}\label{la_n}
\lambda_{n,j}=(2n\pi)^{2}+o(1),\qquad j=1,2,\ \ \text{as}\ \,n\rightarrow
\infty.
\end{equation}

For large indices $n$ the numbers $\lambda_{n,1}$ and
$\lambda_{n,2}$ represent simple eigenvalues, provided that they
do not coincide. Otherwise the multiplicities of these eigenvalues
equal $2$ and the corresponding root subspaces consist either of
two eigenfunctions or of Jordan chains which are formed by one
eigen and one associated functions. Of course for the first
indices the structure of the root subspaces can be more
complicated. The number of such indices is finite, and for
convenience without loss of generality  we assume that the
multiplicity of all eigenvalues $\leqslant2$.

To simplify the notations we shall omit the index $j$ and
enumerate the eigenvalues in the following way:
$\lambda_{n,1}=:\lambda_{n},\ \,\lambda _{n,2}=:\lambda_{-n}$. We
remark that there is a one-to-one correspondence between the
eigenvalues (counting with multiplicities) and integers which
preserves asymptotic (3). In other words: there exists a number $n_{0}%
\geqslant1$, such that for all $n\geqslant n_{0}$ each disk $|\lambda-(2\pi
n)^{2}|\leqslant1$ contains exactly 2 eigenvalues $\lambda_{n}$
and $\lambda_{-n}$ (counting with multiplicities), and the number
of all the remaining eigenvalues equals $2(n_{0}-1)+1$. This
statement is obvious for the potential $q(x)=0$. In a general case
it can  easily be  proved in a standard way by using Rouche's
theorem (we omit the proof of this fact, since it is used only to
simplify the notations).

 For large multiple eigenvalues we have $\lambda
_{n}=\lambda_{-n}$, the remaining ones can be enumerated in such a
way that the same relation holds. This follows from the previous
remark and the agreement that all the eigenvalues have
multiplicity $\leqslant 2$. Denote by $\varphi_{\pm n}$ the
eigenfunction corresponding to simple eigenvalues $\lambda_{\pm
n}$. In the case when the eigenvalue $\lambda_{n}=\lambda_{-n}$
has geometric multiplicity 2 we choose the pair of normalized
eigenfunctions $\varphi_{n},\varphi_{-n}$ such that they are
mutually orthogonal. In the case when the only eigenfunction
$\varphi_{n}$ corresponds to the multiple eigenvalue
$\lambda_{n}=\lambda _{-n}$ we assume that
$\Vert\varphi_{n}\Vert=1$ and choose the associate function
$\varphi_{-n}=:\psi_{n}$ to be orthogonal to $\varphi_{n}$ (it is
uniquely defined by this condition). The obtained system
$\{\varphi _{n}\}_{-\infty}^{\infty}$ of the eigen- and associated
functions (EAF) will be called
\textit{normal}.
The set of positive indices $n\in\mathbb{N}$ which correspond to
the simple (multiple) eigenvalues $\lambda_{n}$ we denote by
$\mathbb{N}_{1}$ $\,$($\mathbb{N}_{2}$). The subset of
$n\in\mathbb{N}_{2}$ which enumerate the Jordan chains
$\varphi_{n},\psi_{n}$ we denote by $\mathbb{N}_{3}$.

Define the Fourier coefficients of the eigenfunctions
\begin{equation}\label{Four}
u_{n}=(\varphi_{n},e^{2\pi i|n|x}),\quad v_{n}=(\varphi_{n},e^{-2\pi
i|n|x}),\quad n=0,\pm1,\pm2,\dots.
\end{equation}

Integrating the first term in the equality
\[
(-\varphi_{n}^{\prime\prime}+q(x)\varphi_{n}-\lambda_{n}\varphi_{n},\,e^{2\pi
ikx})=0
\]
by parts we get
\begin{equation}\label{simple}
\lbrack\lambda_{n}-(2\pi k)^{2}](\varphi_{n}(x),e^{2\pi ikx})=(q(x)\varphi
_{n}(x),e^{2\pi ikx}),
\end{equation}
for all $k\in\mathbb{Z.}$ Similarly, if the Jordan chain
$\varphi_{n},\psi_{n}$, corresponds to the eigenvalue
$\lambda_{n}$ then the equalities
\begin{gather}
\left(  \left[  -\psi_{n}^{\prime\prime}+q(x)\psi_{n}-\lambda_{n}\psi
_{n}-\varphi_{n}\right]  ,\ e^{2\pi ikx}\right)  =0,\nonumber\\
\label{assimple}\left[  \lambda_{n}-(2\pi k)^{2}\right]  (\psi_{n}(x),\ e^{2\pi ikx}%
)=(q(x)\psi_{n}(x),\ e^{2\pi ikx})-(\varphi_{n}(x),\ e^{2\pi ikx})
\end{gather}
are valid. It follows from asymptotics (3) that for large
$|n|>n_{0},$ and $n\neq\pm k$ the estimates
\begin{equation}\label{ineq}
\left\vert \lambda_{n}-(2\pi k)^{2}\right\vert \geq|(n-k)(n+k)|
\end{equation}
hold. Now using \eqref{simple}-\eqref{ineq} we shall
obtain the asymptotics for the eigen- and associated
functions. We do not consider the statement of the
subsequent lemma as a new one (see related results in
\cite[Ch.1]{Na}, \cite{VD}, \cite[Proposition
4]{DM2}, \cite[Proposition 11]{DM1}, for example).
However, we could not find a proper reference to all
assertions of the lemma and decided to present here a
short proof.

\begin{lemma}
The eigenfunctions  $\varphi_{n}$ admit a representation of the
form
\begin{equation}\label{EF}
\varphi_{n}(x)=u_{n}e^{2\pi inx}+v_{n}e^{-2\pi inx}+\varphi_{n}^{0}(x),
\end{equation}
where $u_n, v_n$  are the Fourier coefficients defined by
\eqref{Four}, the functions $\varphi^0_n(x) \perp e^{\pm2\pi inx}$, and
obey the estimates
$$
\sup\limits_{x\in\lbrack
0,1]} |\varphi_{n}^{0}(x)|=O\left(  n^{-1}\ln n\right),
\quad \Vert
\varphi_{n}^{0}(x)\Vert=O(n^{-1}), \ \, n\in\mathbb{Z},
$$
provided that  $\|\varphi_n\| =1$.

Analogously, if (for $n\in\mathbb N_3$) a function $\psi_n(x)$ is
associate with $\varphi_n(x)$, then a representation of the form
\begin{equation}\label{AF}
\psi_{n}(x)=\tilde{u}_{n}e^{2\pi inx}+\tilde{v}_{n}e^{-2\pi inx}+
\psi_n^0(x),
\end{equation}
is valid where $\tilde{u}_{n}=(\psi_{n},e^{2\pi
inx}),\quad\tilde{v}_{n}=(\psi _{n},e^{-2\pi inx}),
\quad \psi^0_n
\perp e^{\pm2\pi inx}$, and the estimates
$$
\sup\limits_{x\in\lbrack
0,1]} |\psi_{n}^{0}(x)|=O\left(  n^{-1}\ln n\right)(1+\|\psi_n\|),
\quad \Vert
\psi_{n}^{0}(x)\Vert=O(n^{-1})(1+\|\psi_n\|), \ \, n\in\mathbb{Z}.
$$
hold.
\end{lemma}

\begin{proof}
Since $q\varphi_{n}\in L_{1}[0,1],$ we have
\[
\lim_{k\rightarrow\infty}(q(x)\varphi_{n}(x),e^{2\pi ikx})=0.
\]
Therefore there exist $C_{n}$ and $k_{0}$ such that%
\[
\max_{k}\mid(q(x)\varphi_{n}(x),e^{2\pi ikx})\mid=\mid(q(x)\varphi
_{n}(x),e^{2\pi ik_{0}x})\mid=C_{n}.
\]
Using  relations  \eqref{simple} and \eqref{ineq} we find
\begin{multline*}
|\varphi_{n}^{0}(x)|\leq\sum_{k\neq\pm n}|(\varphi_{n}(x),e^{2\pi ikx}%
)|\leq\sum_{k\neq\pm n}\left\vert \frac{(q(x)\varphi_{n}(x),e^{2\pi ikx}%
)}{\lambda_{n}-(2\pi k)^{2}}\right\vert =\\
=O(C_{n})\sum_{k\neq\pm n}|n-k|^{-1}\ |n+k|^{-1}=O\left(  C_{n}\frac{\ln n}%
{n}\right).
\end{multline*}
Similarly we get the estimate
\[
\Vert\varphi_{n}^{0}(x)\Vert^{2}\leq\sum_{k\neq\pm n}\left\vert \frac
{(q(x)\varphi_{n}(x),e^{2\pi ikx})}{\lambda_{n}-(2\pi k)^{2}}\right\vert
^{2}=O(C_{n})\sum_{k\neq\pm n}|n-k|^{-2}\ |n+k|^{-2}=O\left(  \frac{C_{n}%
}{n^{2}}\right).
\]
Further, taking into account the inequalities $|u_n|^2 + |v_n|^2
\leqslant 1$ and using representation \eqref{EF}, we find
\begin{multline*}
\mid C_{n}\mid = \mid(q(x)\varphi_{n}(x),e^{2\pi ik_{0}x})\mid=\mid(q(x)(u_{n}%
e^{2\pi inx}+v_{n}e^{-2\pi inx}+\varphi_{n}^{0}(x)),e^{2\pi
ik_{0}x})\mid\leqslant
\\
\left( 2+O\left(  C_{n}n^{-1}\ln n\right)  \right)  \int_{0}^{1}%
\mid(q(x)\mid dx.
\end{multline*}
This relation implies  $C_n =O(1)$. Hence, we have proved the
estimates for the functions $\varphi_n^0$. The estimates for the
functions  $\psi_n^0$ in representation \eqref{AF} can be obtained
similarly taking into account equalities
\eqref{assimple}.
\end{proof}

Recall that a system $\{f_{k}\}$ of elements in a Hilbert space
$H$ is said to be \textit{a Riesz basis }of this space if it is
equivalent to some (and then to all) orthonormal basis
$\{e_{k}\}$. The equivalence means that there is a bounded and
boundedly invertible operator $A$ such that $Ae_{k}=f_{k}$. A
basis $\{f_{k}\}$ is called \textit{unconditional}, if it remains
to be a basis after all rearrangements of its elements. It is
known \cite[Ch.6]{GK}, that a basis is unconditional if and only
if it is equivalent to an orthogonal one. One says that
$\{f_{k}\}$ is
\textit{a Bessel system}, if the series $\sum|(f,f_{k})|^{2}$
converges for all elements $f$ in $H.$ The notation\ $a_{n}\asymp
b_{n}$ means further that there exist positive constants $c_{1},$
$c_{2}$ such that $c_{1}|a_n|<b_{n}<c_{2}|a_n| $.

The following statement plays an important role in the proof of
main results.

\begin{theorem}
 The following assertions are equivalent:

\begin{itemize}
\item [(i)] A normal system of EAF of the operator $L_{0}$
forms a Riesz basis in the space $L_2[0,1]$;

\item [(ii)] The set of indices $\mathbb{N}_{3}$ (i.e. the number of Jordan
chains) is finite and the relations
\begin{gather}\label{n+}
u_{n}\asymp v_{n} \qquad n\in\mathbb{N}_{1},\\
\label{n-}
u_{-n}\asymp v_{-n}\quad n\in\mathbb{N}_{1},
\end{gather}
hold, where $u_{n},v_{n}$ are the Fourier coefficients defined in (4);

\item [(iii)] The set of indices $\mathbb{N}_{3}$ is finite and either relation
(10) or (11) holds.
\end{itemize}
\end{theorem}

\begin{proof}
 \textit{Step 1.} Assume that the number of Jordan chains $\varphi_{n}%
,\psi_{n},n\in\mathbb{N}_{3},$ is infinite. According to the definition of a
normal system we have $\Vert\varphi_{n}\Vert=1$. Let us show that $\Vert
\psi_{n}\Vert\rightarrow\infty$ as $n\rightarrow\infty$. This will imply
 that the condition for the set $\mathbb{N}_{3}$ to be finite is
necessary for the Riesz basis property of the EAF system (since
the condition for the norms of elements to be bounded from above
and below is necessary for the Riesz basis).

Let us substitute expression (9) for the function $\psi_{n}$ into (6) putting
$k=n$ and $k=-n$. Then we get
\begin{gather*}
\lbrack\lambda_{n}-(2\pi n)^{2}]\tilde{u}_{n}-q_{2n}\tilde{v}_{n}%
+u_{n}=O(n^{-1}\ln n)(1+\Vert\psi_{n}\Vert),\\
\lbrack\lambda_{n}-(2\pi n)^{2}]\tilde{v}_{n}-q_{-2n}\tilde{u}_{n}%
+v_{n}=O(n^{-1}\ln n)(1+\Vert\psi_{n}\Vert).
\end{gather*}
where $q_{k}=(q(x),e^{2\pi ikx})$ are the Fourier coefficients of the
potential $q\in L_{1}[0,1]$. We have $q_{2n},q_{-2n}\rightarrow0$ and
according to (3), $\lambda_{n}-(2\pi n)^{2}\rightarrow0$ as $n\rightarrow
\infty$. Suppose $\Vert\psi_{n}\Vert<C$ for some subsequence of indices
$n$. Since the third function in the representation (9) is
orthogonal to the first two ones, we find $|\tilde
{u}_{n}|^{2}+|\tilde{v}_{n}|^{2}<C^{2}$. Then it follows from the obtained
equalities $|u_{n}|+|v_{n}|\rightarrow0$. This is a contradiction to the
condition $\Vert\varphi_{n}\Vert=1$. Consequently, $\Vert\psi_{n}%
\Vert\rightarrow\infty.$

\textit{Step 2.} Let the system $\{\varphi_{n}\}_{-\infty}^{\infty}$ forms a
Riesz basis. Then the number of Jordan chains is finite. For
simplicity and without loss of generality we may assume that
$\mathbb N_3=\emptyset$. Remark that the functions
$\overline{\varphi}_{n}(x),(n\notin\mathbb{N}_{3})$ are the
eigenfunctions of the adjoint operator $L_{0}^{\ast}$ (it is
generated by the same differential expression with the potential
$\overline{q}$ and periodic boundary conditions.) It follows from
the biorthogonality relations (see \cite[ Ch. 1]{Ma}, for example)
that for all $n\in\mathbb{N}_{1}$, the relations
$(\varphi_{n},\,\overline{\varphi}_{k})=\alpha_{n}\delta_{nk},$
hold. Here $\alpha_{n}\neq0$, and $\delta_{nk}$ is the Kroneker
symbol. Using representation \eqref{EF} and the estimate for the
function $\varphi_n^0$ in $L_2$-norm, we find
\begin{equation}\label{al_n}
\alpha_{n}=(\varphi_{n},\,\overline{\varphi}_{n})=2u_{n}v_{n}+O(n^{-2}),
\quad n\in\mathbb N_1.
\end{equation}
The eigenfunctions $\varphi_{n}$ and $\varphi_{-n}$ corresponding
to the indices
$n\in\mathbb{N}_{2}\setminus\mathbb{N}_{3}=\mathbb{N}_{2}$ are
chosen mutually orthogonal. Then, since the function $\varphi_{n}$
has representation \eqref{EF}, we have
\[
\varphi_{-n}(x)=\overline{v}_{n}e^{2\pi inx}-\overline{u}_{n}e^{-2\pi
inx}+O(n^{-1})
\]
up to a scalar multiple of modulus 1. Consequently, the system
$\{\chi _{n}(x)\}_{-\infty}^{\infty}$ which is biorthogonal to
$\{\varphi _{n}\}_{-\infty}^{\infty}$ admits representation
\[
\chi_{n}=\alpha_{n}^{-1}\overline{\varphi}_{n},\quad\text{if}\ \,n\in
\mathbb{N}_{1},\quad\chi_{\pm n}=\overline{\varphi}_{\mp n}+O(n^{-1}%
),\quad\text{if}\ \,n\in\mathbb{N}_{2}.
\]
The so-called uniform minimality condition $\Vert\varphi_{n}\Vert\,\Vert
\chi_{n}\Vert\leq const$ is necessary for the basis property. Therefore,
$|\alpha_{n}|^{-1}\Vert\varphi_{n}\Vert^{2}=|\alpha_{n}|^{-1}\leq
\text{const},n\in\mathbb{N}_{1}.$ Using this estimate, equalities
\eqref{al_n} and the relations
$|u_{n}|^{2}+|v_{n}|^{2}=1+O(n^{-1})$ we get $u_{n}\asymp v_{n},n\in
\mathbb{N}_{1}$. Thus, we have proved the implication (i)
$\Rightarrow $(iii).

\textit{Step 3.}
Let either condition \eqref{n+} or \eqref{n-} hold. Assume, for
example, that the first one  $u_n\asymp v_n \asymp 1$ is true.  It
follows from the biorthogonality relations that
\[
0=(\varphi_{n},\,\overline{\varphi}_{-n})=u_{n}v_{-n}+v_{n}u_{-n}+O(n^{-1}).
\]
Since  $|u_{-n}| +|v_{-n}| \asymp 1$, the last relation holds only
if the condition
\eqref{n-} is valid. Hence, we have proved the equivalence
 $(ii)\Leftrightarrow (iii).$

The proof of the implication $(ii)
\Rightarrow (i)$ can be found in the papers
\cite{DV} and \cite{M}. We remark only that it can be readily
obtained independently: condition (ii) together with obtained
relations \eqref{al_n} imply that the angles between the
normalized functions
 $\varphi_n$ and $\varphi_{-n}$ are uniformly positive. Indeed, if
 the angles tend to zero for some indices $n\to\infty$
 then $\varphi_n=\varphi_{-n} +o(1)$ for these indices
 and relations $\alpha_n = (\varphi_n, \overline\varphi_n),
 (\varphi_n,\overline\varphi_{-n})=0$ imply $\alpha_n \to 0$.
 This contradicts to \eqref{al_n}.  Then the
implication $(ii)\Rightarrow (i)$ follows from theorem
\cite{Sh1} on the Riesz basis from 2-dimensional subspaces.
This completes the proof of Theorem.
\end{proof}

\begin{remark} In the case when the number
of multiple eigenvalues is finite the equivalence of the first two
statements was used in papers \cite{DV} and \cite{M}. Here we used
the arguments from \cite{DV}. We remark that in the case of
multiple eigenvalues the assertion of Theorem~1 from \cite{M} has
to be corrected.
\end{remark}
\medskip

Let us formulate another statement which gives criteria for unconditional
basis property. Its proof repeats the arguments of the previous theorem and is
omitted here.

\noindent\textbf{Theorem $1^{\prime}$.} \emph{The system $\{\varphi_{n}\}$ of EAF of
the operator $L_{0}$ forms an unconditional basis \ if \ and only
\ if \  either condition (10) or (11) }holds.

\section{Main results}

The main results will be formulated later, first we prove several
lemmata. The main tool in our approach to solve the basis problem
is the relation \eqref{basic} obtained in Lemma~3. However, the
most essential part of the paper is to encode in terms of the
potential $q$ the behavior of numbers in the sequences
$B_{m}(\lambda_{n})$ and $B_{m}^{\prime}(\lambda_{n})$
participating in \eqref{basic}.

The first of the subsequent lemmata is obtained in the paper
\cite{DV}. Here we only remark that its proof is based on formula
\eqref{simple}. One has to put $k=n$ in this formula to obtain
\eqref{2.1}
and $k=-n$ to obtain \eqref{2.2}. In the first case the right
hand-side in formula \eqref{simple} has to be replaced by the
series (we use the assumption $q_{0}=0$)
\[
\sum_{n_{1}\neq0}q_{n_{1}}(\varphi_{n}(x),e^{2\pi i(n-n_{1})x}),
\]
where $q_{n}$ are the Fourier coefficients of the potential $q$.
It can be shown that the series converges and represent the right
hand-side. Then, we do not change the terms with indices
$n_{1}=2n$; all the other terms for $n_{1}\neq2n$ we replace
according to \eqref{simple} by the expressions
\[
q_{n_{1}}\frac{(q(x)\varphi_{n}(x),e^{2\pi i(n-n_{1})x})}{\lambda_{n}%
-(2\pi(n-n_{1}))^{2}}.
\]
The same procedure can be applied to the numerator of the last
expression, and then one can proceed in a similar way.

\begin{lemma}
Let $q_{n}=(q(x),\,e^{2\pi inx})$ and $q_{0}=0$. The following relations are
valid ($n\geq1$)
\begin{equation}\label{2.1}
\lbrack\lambda_{n}-(2\pi n)^{2}-A_{m}(\lambda_{n})]\ (\varphi_{n}(x),\,e^{2\pi
inx})-(q_{2n}+B_{m}(\lambda_{n}))\ (\varphi_{n},\,e^{-2\pi inx})=R_{m},
\end{equation}
where
\begin{gather*}
A_{m}(\lambda_{n})=\sum_{k=1}^{m}a_{k}(\lambda_{)},\qquad B_{m}(\lambda
_{n})=\sum_{k=1}^{m}b_{k}(\lambda_{n}),\\
a_{k}(\lambda_{n})=\sum_{n_{1},n_{2},...,n_{k}}\frac{q_{n_{1}}q_{n_{2}}\cdots
q_{n_{k}}q_{-n_{1}-n_{2}-...-n_{k}}}{[\lambda_{n}-(2\pi(n-n_{1}))^{2}%
]\cdots\lbrack\lambda_{n}-(2\pi(n-n_{1}-\dots-n_{k}))^{2}]},\\
b_{k}(\lambda_{n})=\sum_{n_{1},n_{2},\dots,n_{k}}\frac{q_{n_{1}}q_{n_{2}%
}\cdots q_{n_{k}}q_{2n-n_{1}-n_{2}-\dots-n_{k}}}{[\lambda_{n}-(2\pi
(n-n_{1}))^{2}]\cdots\lbrack\lambda_{n}-(2\pi(n-n_{1}-\dots-n_{k}))^{2}]},\\
\ \ \ R_{m}=\sum_{n_{1},n_{2},\dots,n_{m+1}}\frac{q_{n_{1}}q_{n_{2}}\cdots
q_{n_{m+1}}(q(x)\varphi_{n}(x),\,e^{2\pi i(n-n_{1}-\dots-n_{m+1})})}%
{[\lambda_{n}-(2\pi(n-n_{1}))^{2}]\cdots\lbrack\lambda_{n}-(2\pi(n-n_{1}%
-\dots-n_{m+1}))^{2}]}.
\end{gather*}
The summation in these formulae is taken over the indices
$n_{p}\neq 0,\sum_{1}^{p}n_{j}\neq0,\sum_{1}^{p}n_{j}\neq2n$ for
$1\leq p\leq m+1$.

The relations
\begin{equation}\label{2.2}
\lbrack\lambda_{n}-(2\pi n)^{2}-A_{m}^{\prime}(\lambda_{n})]\ (\varphi
_{n}(x),\,e^{-2\pi inx})-(q_{-2n}+B_{m}^{\prime}(\lambda_{n}))\ (\varphi
_{n},\,e^{2\pi inx})=R_{m}^{\prime}%
\end{equation}
are valid, too. Here
\begin{gather*}
A_{m}^{\prime}(\lambda_{n})=\sum_{k=1}^{m}a_{k}^{\prime}(\lambda_{)},\qquad
B_{m}^{\prime}(\lambda_{n})=\sum_{k=1}^{m}b_{k}^{\prime}(\lambda_{n}),\\
a_{k}^{\prime}(\lambda_{n})=\sum_{n_{1},n_{2},\dots,n_{k}}\frac{q_{n_{1}%
}q_{n_{2}}\cdots q_{n_{k}}q_{-n_{1}-n_{2}-\dots-n_{k}}}{[\lambda_{n}%
-(2\pi(n+n_{1}))^{2}]\cdots\lbrack\lambda_{n}-(2\pi(n+n_{1}+\dots+n_{k}%
))^{2}]},\\
b_{k}^{\prime}(\lambda_{n})=\sum_{n_{1},n_{2},\dots,n_{k}}\frac{q_{n_{1}%
}q_{n_{2}}\cdots q_{n_{k}}q_{-2n-n_{1}-n_{2}-\dots-n_{k}}}{[\lambda_{n}%
-(2\pi(n+n_{1}))^{2}]\dots\lbrack\lambda_{n}-(2\pi(n+n_{1}+\dots+n_{k}))^{2}%
]},\ \\
\ \ \ R_{m}^{\prime}=\sum_{n_{1},n_{2},\dots,n_{m+1}}\frac{q_{n_{1}}q_{n_{2}%
}\cdots q_{n_{m+1}}(q(x)\varphi_{n}(x),\,e^{2\pi i(n+n_{1}+\dots+n_{m+1})}%
)}{[\lambda_{n}-(2\pi(n+n_{1}))^{2}]\cdots\lbrack\lambda_{n}-(2\pi
(n+n_{1}-\dots+n_{m+1}))^{2}]},
\end{gather*}
and the summation is taken over the indices $n_{p}\neq0,\ \sum_{1}^{p}%
n_{j}\neq0,\ \sum_{1}^{p}n_{j}\neq-2n,\ 1\leq p\leq m+1$.

The terms in the above formulae admit the estimates
\[
a_{k}(\lambda_{n})=O\left(  \frac{\ln^{k} \,
n}{n^{k}}\right)  ,\quad b_{k}(\lambda_{n})=O\left(
\frac{\ln^{k} \, n}{n^{k}}\right)  , \quad
R_{m}=O\left(  \frac{\ln^{m+1}\, n}{ n^{m+1}}\right)
,
\]
and the same estimates admit the terms $a^{\prime}_{k}(\lambda_{n}%
),\ b^{\prime}_{k}(\lambda_{n})$ and $R^{\prime}_{m}$.

The above formulae remain valid for the negative indices $-n\notin
\mathbb{N}_{3}$ if \ $\ln\,n$ is replaced by\ $\ln\,|n|$.
They remain also valid
if the eigenvalues $\lambda_{n}$ are replaced by $\lambda_{-n}$.
\end{lemma}

\begin{lemma}
Let $u_{n},v_{n},\ n\geq1,$ be the Fourier coefficients of the
eigenfunctions defined according to formulae \eqref{Four}. Let
$m\geq0$ be arbitrary integer. Then the relation
\begin{equation}\label{basic}
(q_{2n}+B_{m}(\lambda_{n}))v_{n}^{2}-(q_{-2n}+B_{m}^{\prime}(\lambda
_{n}))u_{n}^{2}=O\left(  \frac{\ln^{m+1}\,n}{n^{m+1}}\right)
\end{equation}
is valid.
\end{lemma}

\begin{proof}
Let us rewrite relations \eqref{2.1} and \eqref{2.2} in the form
\begin{gather*}
(\lambda_{n}-(2\pi
n)^{2}-A_{m}(\lambda_{n}))u_{n}-(q_{2n}+B_{m}(\lambda
_{n}))v_{n}=O\left(  \frac{\ln^{m+1}\,n}{n^{m+1}}\right)  ,\\
(\lambda_{n}-(2\pi n)^{2}-A_{m}^{\prime}(\lambda_{n}))v_{n}-(q_{-2n}%
+B_{m}^{\prime}(\lambda_{n}))u_{n}=O\left(  \frac{\ln^{m+1}\,n}{n^{m+1}%
}\right)  .
\end{gather*}
We shall prove the equalities $A_{m}(\lambda_{n})=A_{m}^{\prime}(\lambda_{n}%
)$. Then, relation \eqref{basic} is obtained from the last two
ones: multiply the first and the second relation by $v_{n}$ and
$u_{n}$, respectively, and take the difference.

It is sufficient to prove that $a_{k}(\lambda_{n})=a_{k}^{\prime}(\lambda
_{n})$. Let us make the substitution
\[
-n_{1}-n_{2}-\dots-n_{k}=j_{1},\ \,n_{2}=j_{k},\ \,n_{3}=j_{k-1}%
,\,\dots,\,n_{k}=j_{2},
\]
in the formula for the expression $a_{k}^{\prime}(\lambda_{n}).$
Then the inequalities for the forbidden indices
$n_{p}\neq0,\sum_{1}^{p}n_{s}\neq0,-2n$ for $1\leq p\leq k$ in the
formula for $a_{k}^{\prime}$ take the form
$j_{p}\neq0,\sum_{1}^{p}j_{s}\neq0,2n$ for $1\leq p\leq k$, and it
will coincide with the formula for $a_{k}.$ The lemma is proved.
\end{proof}

\begin{lemma}
 Let
 for some $m\geq0$ one of the following relations hold
\begin{equation}\label{o}
(q_{2n}+B_{m}(\lambda_{n}))^{-1}=o\left(  \frac{n^{m+1}}
{\ln^{m+1}\,n}\right),
\quad(q_{-2n}+B_{m}^{\prime}(\lambda_{n}))^{-1}=o\left(  \frac{n^{m+1}}%
{\ln^{m+1}\,n}\right),
\end{equation}
as $ n\to+\infty$.
 Then the number of the eigenvalues with geometric multiplicity 2
is finite (i.e. the set $\mathbb{N}_{2}\setminus\mathbb{N}_{3}$ is
finite).
\end{lemma}

\begin{proof}
Let
\[
\varphi_{n}(x)=u_{n}e^{2\pi inx}+v_{n}e^{-2\pi inx}+\varphi_{n}^{0}%
(x),\quad\varphi_{-n}(x)=u_{-n}e^{2\pi inx}+v_{-n}e^{-2\pi inx}+\varphi
_{-n}^{0}(x),
\]
be the eigenfunctions corresponding to an eigenvalue $\lambda_{n}=\lambda
_{-n}$. By Lemma~1 we have $\Vert\varphi_{n}^{0}\Vert+\Vert\varphi_{-n}%
^{0}\Vert=O(n^{-1})$. It is assumed that $\varphi_{n}$ and
$\varphi_{-n}$ are chosen  orthogonal. Therefore, the determinant
which is compiled from the rows
$\{u_{n},v_{n}\},\{u_{-n},v_{-n}\}$ is not equal to zero. Hence,
there are linear combinations of these functions having the same
representations with $u_{n}=0,\ \,v_{n}=1+O(n^{-1}),\
\,u_{-n}=1+O(n^{-1}),\ \,v_{-n}=0$. The relation \eqref{basic}
takes for these functions the form
\[
q_{2n}+B_{m}(\lambda_{n})=O\left(  \frac{\ln^{m+1}\,n}{n^{m+1}}\right)  ,\quad
q_{-2n}+B_{m}^{\prime}(\lambda_{n})=O\left(  \frac{\ln^{m+1}\,n}{n^{m+1}%
}\right)  \quad\text{as}\ n\to\infty.
\]
Consequently, no of two relations in \eqref{o} can be valid if the
number of the eigenvalues of geometric multiplicity 2 is infinite.
The lemma is proved.
\end{proof}

\begin{lemma} 
Let the both estimates in \eqref{o} are valid for some $m\geq0$,
and in addition the estimates
\begin{equation}\label{above}
(q_{2n}+B_{m}(\lambda_{n}))= O\left(\frac{\ln^{m+1}\,
n}{n^{m-3}}\right), \quad
(q_{-2n}+B_{m}^{\prime}(\lambda_{n}))^{-1}=O\left(
\frac{\ln^{m+1}\, n}{n^{m-3}}\right)
\end{equation}
hold, as  $n\to+\infty$ (for $m\leq3$ they hold automatically).
Then the number of associated functions (i.e. the number of
indices $\mathbb N_3$ is finite).
\end{lemma}

\begin{proof}
Let $\varphi_{n},\psi_{n}$ be a Jordan chain corresponding to the
eigenvalue $\lambda_{n}$. Then
$\overline{\varphi}_{n},\overline{\psi}_{n}$ is the Jordan chain
of the adjoint operator $L_{0}^{\ast}$ corresponding to the
eigenvalue $\overline{\lambda}_{n}$. The biorthogonality relations
(see \cite[Ch.1]{Ma}, for example) give
$(\varphi_{n},\overline{\varphi}_{n})=0$. Using equality
\eqref{al_n}, we get
\[
u_{n}v_{n}=O(n^{-2}).
\]
The condition for the norms
$\|u_{n}\|^{2}+\|v_{n}\|^{2}=1+O(n^{-1})$ implies that either
$\|v_{n}\|>1/2$ or $\|u_{n}\|>1/2$ for all sufficiently large $n$.
Let, for example, the first inequality holds for  infinite set of
indices $n$. Then the estimate $u_{n}=O(n^{-2})$ holds for these
indices, and relation
\eqref{basic} together with the  second estimate in \eqref{above} give
\[
(q_{2n}+B_{m}(\lambda_{n}))=
(q_{-2n}+B_{m}^{\prime}(\lambda_{n}))O(n^{-4} )+O\left(
\frac{\ln^{m+1}\,n}{n^{m+1}}\right) =O\left(  \frac{\ln^{m+1}
\,n}{n^{m+1}}\right) .
\]
This contradicts to condition \eqref{o}. The case when
$\|u_{n}\|>1/2$ for infinite set of indices $n$ can be treated
analogously.  The lemma is proved.
\end{proof}

\begin{lemma}
Let $p\geq0$ be an arbitrary integer,
\begin{equation}\label{condition}
q(x)\in W_{1}^{p}[0,1],\quad\text{ and}
\ \, q^{(l)}(0)=q^{(l)}(1)\quad\text{for
all}\ \,0\leq l\leq s-1,\text{ with some }s\leq p.
\end{equation}
Define the functions
\[
Q(x)=\int_{0}^{x}q(t)\,dt,\quad S(x)=Q^{2}(x),
\]
and denote by $Q_{k}=(q(x),\,e^{2\pi ikx}),\
\, S_{k}=(q(x),\,e^{2\pi ikx})$ the Fourier coefficients of these
functions.

The following relations are valid:
\begin{gather}\label{19}
b_{1}(\lambda_{n})=-S_{2n}+2Q_{0}Q_{2n}+o\left(  n^{-s-2}\right)
,\quad b_{1}^{\prime}(\lambda_{n})=-S_{-2n}+2Q_{0}Q_{-2n}+o\left(
n^{-s-2}\right) ,\\
\label{20}
b_{2}(\lambda_{n,j})=o\left(
n^{-s-2}\right)  ,\quad b_{2}^{^{\prime} }(\lambda_{n})=o\left(
n^{-s-2}\right)  ,
\\ \label{21}
b_{k}(\lambda_{n})=o\left(  \frac{\ln^{k}\,n}{n^{k+s}}\right)
,\quad b_{k}^{^{\prime}}(\lambda_{n})=o\left(
\frac{\ln^{k}\,n}{n^{k+s}}\right)
\quad\text{for all} \  k\geq 3,4,\dots.
\end{gather}
\end{lemma}

\begin{proof} {\sl Step 1.}
Integrating the equalities $q_{n}=(q(x),e^{2\pi inx})$ by parts
 and using the assumptions of the lemma, we get
$q_{n}=o(n^{-s})$. Further, at least one of the numbers
$n_{1},n_{2},\dots,n_{k}$ and $\pm2n-n_{1}-n_{2}-\dots-n_{k}$ has
modulus greater then $n/k$. Therefore,
\[
q_{n_{1}}q_{n_{2}}\cdots
q_{n_{k}}q_{\pm2n-n_{1}-n_{2}-\dots-n_{k}} =o(n^{-s}).
\]
The denominators in the expressions for $b_{k}(\lambda_{n})$ and
$b_{k}^{\prime}(\lambda_{n})$ we can estimate from below using the
inequalities \eqref{ineq}. Then we get the series not depending on
$\lambda_{n}$. Its estimate is trivial and we readily get
\eqref{21}.

{\sl Step 2.}
The proof of the first two relations in the assertion of the lemma
requires a sharper analysis. From relation \eqref{la_n} we easily
get the estimate
\[
\sum_{\substack{k\neq0,\pm2n}}^{\infty}\mid\frac{1}{\lambda_{n}-(2\pi(n\mp
k))^{2}}-\frac{1}{(2\pi n)^{2}-(2\pi(n\mp k))^{2}}\mid=o(\frac{1}{n^{2}}).
\]
By virtue of this estimate
\begin{equation}
\label{22}
b_{1}(\lambda_{n})=\widetilde{S}(2n)+o\left(  n^{-s-2}\right)  ,\quad
b_{1}^{^{\prime}}(\lambda_{n})=\widetilde{S}(-2n)+o\left(  n^{-s-2}\right)  ,
\end{equation}
where
\[
\widetilde{S}(\pm2n)=\sum_{k\neq0,\,2n}\frac{q_{k}\,q_{2n-k}}{(2\pi
n)^{2}-(2\pi(n\mp k))^{2}}=\frac{1}{4\pi^{2}}\sum_{k\neq0,\,2n}\frac
{q_{k}\,q_{\pm2n-k}}{k(\pm2n-k)}.
\]
Now, let us substitute in the equality $S_{\pm2n}=(Q^{2}(x),e^{\pm2\pi
i(2n)x})$ the Fourier series
\[
Q(x)=\sum_{k\in\mathbb{Z}}Q_{k}e^{2\pi ikx}
\]
and view in mind that $Q_{k}=(2\pi ik)^{-1}q_{k}$ for $k\neq0.$ Then we get
\begin{equation*}
\widetilde{S}(2n)=-S_{2n}+2Q_{0}Q_{2n},\quad \widetilde{S}(-2n)=-S_{-2n}
+2Q_{0}Q_{-2n}.
\end{equation*}
These estimates together with relations \eqref{22} give relations
\eqref{19}.

{\sl Step 3.}
Let us prove \eqref{20}. Arguing as in the proof of equalities
\eqref{19}, we find
\begin{equation*}
b_{2}(\lambda_{n})=\frac{1}{(2\pi)^{4}}C(2n)+o\left(  \frac{1}{n^{s+3}%
}\right)  ,
\end{equation*}
where
\[
C(2n)=\sum\frac{q_{n_{1}}q_{n_{2}}q_{2n-n_{1}-n_{2}}}{n_{1}(2n-n_{1}%
)(n_{1}+n_{2})(2n-n_{1}-n_{2}))} .
\]
Here and further we imply that the summation is taken over all
indices $n_{1},n_{2}$ which do not annulate the denominator.

Making the substitutions $k_{1}=n_{1},$ $k_{2}=2n-n_{1}-n_{2},$ we
simplify this expression and obtain
\[
C(2n)=\sum\frac{q_{n_{1}}q_{n_{2}}q_{2n-n_{1}-n_{2}}}{n_{1}(2n-n_{1}%
)n_{2}(2n-n_{2}))}.
\]

Using the equality
\[
\frac{1}{k(2n-k)}=\frac{1}{2n}\left(\frac{1}{k}+\frac
{1}{(2n-k)}\right)
\]
we find that
\begin{equation}
C(2n)=\frac{1}{4n^{2}}(I_{1}+2I_{2}+I_{3})
\end{equation}
where%
\[
I_{1}:=\sum\frac{q_{n_{1}}q_{n_{2}}q_{2n-n_{1}-n_{2}}}{n_{1}n_{2}},\text{
}I_{2}:=\sum\frac{q_{n_{1}}q_{n_{2}}q_{2n-n_{1}-n_{2}}}{n_{2}(2n-n_{1}%
)},\text{ }I_{3}:=\sum\frac{q_{n_{1}}q_{n_{2}}q_{2n-n_{1}-n_{2}}}%
{(2n-n_{1})(2n-n_{2})}.
\]
Thus, we have to show that $|I_1|+|I_2|+|I_3| = o(n^{-s})$. To
estimate the sum $I_1$ we remark that
\[
\left((Q(x)-Q_{0})^{2}q(x),e^{2\pi i(2n)x}\right)=-(2\pi)^{-2}I_{1}.
\]
This equality can be readily checked by substituting  in the
left-hand side the Fourier series $Q(x)=Q_{0}+\sum(2\pi
ik)^{-1}q_{k}e^{2\pi ikx}$. Since the function $(Q(x)-Q_0)^2 q(x)$
satisfies condition \eqref{condition}, we obtain
$I_{1}=o(n^{-s})$.

To estimate $I_2$ we remark that
\[
I_{2}=\sum\frac{q_{n_{1}}q_{n_{2}}q_{2n-n_{1}-n_{2}}}{n_{2}(n_{1}+n_{2})}%
=\sum\frac{q_{n_{1}}q_{n_{2}}q_{2n-n_{1}-n_{2}}}{n_{1}(n_{1}+n_{2})}.
\]
Taking the sum of the last expressions we get $2I_2 =I_1$,
therefore $I_2=o(n^{-s})$.

{\sl Step 4}. To estimate  the sum $I_3$ we consider the
functions
\[
G(x, n)=\int_{0}^{x}q(t)e^{-2\pi i(2n)t}dt-q_{2n}x.
\]
Integrating by parts  the expressions for the Fourier coefficients
$ G_k(n): = (G(x,n), e^{2\pi i k x})$  we get
$$
G_k({n}) = \frac 1{2\pi ik } q_{2n+k}, \quad
\text{if} \  k\ne 0,  \quad G(x,n)=G_{0}(n)+\sum_{n_{1}\neq2n}%
\frac{q_{n_{1}}}{2\pi i(n_{1}-2n)}e^{2\pi i(n_{1}-2n)x}.
$$
Using this expansion for the functions $G(x,n)$ we readily find
$$
\int_{0}^{1}(G(x,n)-G_{0}(n))^2 q(x)e^{2\pi i(2n)x}\, dx=\frac{1}{4\pi^{2}}I_{3}.
$$
Now assume that for a function $q(x)$ satisfying condition
\eqref{condition} with $s\geqslant 1$ we have already proved the estimate
\begin{equation} \label{exi}
\int_0^1(G(x,n) -G_0(n)) q^2(x)\, dx =o(n^{-s}).
\end{equation}
Then integration by parts gives (we take into account the estimate $q_{2n} =o(n^{-s})$ and \eqref{exi})
\begin{multline}\label{main}
I_3=4\pi^2 \int_{0}^{1}(G(x,n)-G_{0}(n))^2 q(x)e^{2\pi i(2n)x}\, dx
\\ =\frac {2\pi i}
{n}\int_0^1(G(x,n) -G_0(n))G'(x,n)q(x) e^{2\pi i(2n)x}\, dx\\
 +\frac {\pi i}n\int_0^1
(G(x,n)- G_0(n))^2 q'(x)e^{2\pi i(2n)x}\, dx\\ =\frac {\pi i}n\int_0^1 (G(x,n)- G_0(n))^2 q'(x)e^{2\pi i(2n)x}\, dx
+o(n^{-s-1}).
\end{multline}
It follows from these equations that $I_3 =o(1)$ for  $q$ satisfying condition \eqref{condition} with $s=0$, since any
function $q\in L_1$ can be approximated with  arbitrary accuracy by a smooth periodic function. Then by induction from
\eqref{main} we conclude $I_3 =o(n^{-s})$ that ends the proof.

To prove \eqref{exi} we note that
\begin{multline}\label{I}
G_0(n) =
\int_0^1 \left( \int_0^x q(t)e^{-2\pi i(2n)t}\, dt -x q_{2n}\right)
\, dx=-\int_0^1 x q(x)e^{-2\pi i(2n)x}\, dx +o(n^{-s})\\
 =\frac 1{4\pi i
n }\left[q(1) - \int_0^1(q(x) +x q'(x))e^{-2\pi i(2n)x}\,
dx\right] +o(n^{-s})\\
=\frac {q(1)} {4\pi i n } +\frac {q'(1)}{(4\pi i n)^2} +\dots +
\frac {q^{(s-1)}(1)}{(4\pi i n)^{s-1}}
 +o(n^{(-s)}),
\end{multline}
provided that the first $s-1$ derivatives of $q$ are absolute
continuous functions. Introducing the function $ F(x): =\int_0^x
q^2(x)\, dx$ we also find
\begin{multline} \label{II}
\int_0^xG(x,n)q^2(x)\,dx =-\int_0^1(q(x)e^{-2\pi i(2n)x}
-q_{2n})F(x)\, dx\\
 =\frac 1{4\pi i n }\left[q(1)F(1) -
\int_0^1(q^3(x) +F(x) q'(x))e^{-2\pi i(2n)x}\,
dx\right] +o(n^{(-s)})\\
=F(1)\left[\frac {q(1)} {4\pi i n } +\frac {q'(1)}{(4\pi i n)^2}
+\dots +
\frac {q^{(s-1)}(1)}
{(4\pi i n)^{s-1}}\right]
 +o(n^{(-s)}).
\end{multline}
Equalities \eqref{I} and \eqref{II} give \eqref{exi}. The proof of
the lemma is complete.
\end{proof}

\begin{theorem}
Let one of the relations in \eqref{o} hold. Then a normal EAF
system of the operator $L_{0}$ forms a Riesz basis if and only if
\begin{equation}\label{suf}
q_{2n}+B_{m}(\lambda_{n})\asymp q_{-2n}+B_{m}^{\prime}(\lambda_{n}).
\end{equation}

\end{theorem}

\begin{proof}
Let, for example, the first relation in \eqref{o} hold. Using
relation \eqref{basic}, we obtain
\begin{equation}\label{Th2}
v_{n}^{2}=r_{n}u_{n}^{2}+o(1),\qquad r_{n}:=\frac{q_{-2n}+B_{m}^{\prime
}(\lambda_{n})}{q_{2n}+B_{m}(\lambda_{n})}.
\end{equation}
Assuming the validity of condition \eqref{suf} we obtain that
$r_{n}\asymp
1,\ n\in\mathbb{N}$. Then it follows from \eqref{Th2} that $u_{n}\asymp v_{n}%
\asymp1,\ n\in\mathbb{N}$. This relation implies that the number of Jordan
chains is finite. Indeed (see the proof of Lemma~5), the estimate $u_{n}%
v_{n}=O(n^{-2})$ holds for indices $n\in\mathbb{N}_{3}$ which correspond to
Jordan chains. Therefore, relation $u_{n}\asymp v_{n}\asymp1,\ n\in\mathbb{N}$
may be valid only in the case when the set $\mathbb{N}_{3}$ is finite.
Applying Theorem~1 we find that a normal EAF system forms a Riesz basis.

Conversely, let the system $\{\varphi_{n}\}_{-\infty}^{\infty}$
form a Riesz basis. Then by virtue of Theorem~1 the set of indices
$n\in\mathbb{N}_{3}$ is finite, and it follows from Lemma~4 that
the set of indices $\mathbb{N}_{2}$ is finite, too, i.e all the
eigenvalues are asymptotically simple. In this case Theorem~1
implies that $u_{n}\asymp v_{n}\asymp1,$ and by \eqref{Th2} we
have $r_{n}\asymp1$
\end{proof}

As a corollary we get several results.

\begin{theorem}
Let the potential $q$ obey condition \eqref{condition} and there
is a number $\varepsilon
>0$ such that either the estimate
\begin{equation}\label{assum+}
|q_{2n}-S_{2n}+2Q_{0}Q_{2n}|\geq\varepsilon n^{-s-2}%
\end{equation}
or the estimate
\begin{equation}\label{assum-}
|q_{-2n}-S_{-2n}+2Q_{0}Q_{-2n}|\geq\varepsilon n^{-s-2}%
\end{equation}
hold. Then the condition
\begin{equation}\label{nessuf}
q_{2n}-S_{2n}+2Q_{0}Q_{2n}\asymp q_{-2n}-S_{-2n}+2Q_{0}Q_{-2n}
\end{equation}
is necessary and sufficient for a normal system of EAF to form a Riesz basis.
\end{theorem}

\begin{proof}
It follows from Lemma~6 that for any $m\geq1$
\begin{gather*}
q_{2n}+B_{m}(\lambda_{n})=q_{2n}-S_{2n}+2Q_{0}Q_{2n}+o(n^{-s-2}),\\
q_{-2n}+B_{m}^{\prime}(\lambda_{n})=q_{-2n}-S_{-2n}+2Q_{0}Q_{-2n}+o(n^{-s-2}).
\end{gather*}
Thus, the validity of either relation \eqref{assum+} or
\eqref{assum-} implies the validity of the first or second equality
of \eqref{o}, respectively, with $m\geq s+2$. It is obvious that
relations \eqref{suf} and \eqref{nessuf} are equivalent, provided
that either \eqref{assum+} or \eqref{assum-} holds. It remains to
apply Theorem~2.
\end{proof}

\begin{corollary}
Let the potential $q$ obey condition \eqref{condition} and let one
of the following estimates hold
\[
|q_{2n}|>\varepsilon n^{-s-1},\ \ \text{or}\ \
|q_{-2n}|>\varepsilon n^{-s-1}\ \ \ \forall n\gg1,
\]
with some $\varepsilon>0.$ Then a normal EAF system of the operator $L_{0}$
forms a Riesz basis if and only if $q_{2n}\asymp q_{-2n}$.
\end{corollary}

\begin{proof}
Obviously, $Q,S\in W_{1}^{p+1}[0,1]$ and due to the assumption $q_{0}=0$ we
have
\[
Q^{(l)}(0)=Q^{(l)}(1),\ S^{(l)}(0)=S^{(l)}(1),\quad\text{for }\ \,0\leq l\leq
s.
\]
Then $Q_{\pm2n}=o(n^{-s-1}),S_{\pm2n}=o(n^{-s-1}).$ It remains to apply Theorem~3
\end{proof}

\begin{corollary}
Let $q\in W_{1}^{p}[0,1]$ for some $p\geq1$ and $q^{(s)}(0)\neq q^{(s)}(1)$
for some $s\leq p-1$. Then a normal system of EAF of the operator $L_{0}$ form
a Riesz basis.
\end{corollary}

\begin{proof}
 We may assume that the number $s$ is the smallest one for which
the condition $q^{(s)}(0)\neq q^{(s)}(1)$ holds. This  means that
$q$ obeys condition \eqref{condition}. Integrating by parts we
obtain
\[
q_{n}=\frac{q^{(s)}(1)-q^{(s)}(0)}{(2\pi in)^{s+1}}+o(n^{-s-1}).
\]
This  implies the validity of the estimates from below in
Corollary~1  and $q_{2n}\asymp q_{-2n}.$ Therefore the assertion
follows from Corollary~1.
\end{proof}

The proof of corresponding results for antiperiodic problem can be carried out
in a similar way. Here we only formulate the analogues of Theorem~3,
Corollary~1 and Corollary 2.

\begin{theorem} The following assertions are valid.
\begin{itemize}
\item[(a)]
Let condition \eqref{condition} hold. Assume in addition that
either the estimate
\[
|q_{2n+1}-S_{2n+1}+Q_{0}Q_{2n+1}|\geq n^{-s-2}\qquad\forall\ \,n\gg1
\]
or the estimate
\[
|q_{-2n-1}-S_{-2n-1}+Q_{0}Q_{-2n-1}|\geq n^{-s-2}\qquad\forall\ \,n\gg1,
\]
holds with some $\varepsilon>0$. Then the condition
\[
q_{2n+1}-S_{2n+1}+Q_{0}Q_{2n+1}\asymp q_{-2n-1}-S_{-2n-1}+Q_{0}Q_{-2n-1}%
\]
is necessary and sufficient for a normal EAF system of the
operator $L_1$  to form a Riesz basis.

\item[(b)] Let condition \eqref{condition} hold.
Assume in addition that one of the following  estimates
\[
|q_{2n+1}|\geq\varepsilon n^{-s-1}\quad\forall\ n\gg1\ \,\text{or}\quad
q_{-2n-1}\geq\varepsilon n^{-s-1}\ \,\forall\ n\gg1,
\]
holds with some $\varepsilon>0$. Then the condition
\[
q_{2n+1}\asymp q_{-2n-1}%
\]
is necessary and sufficient for a normal EAF system of the
operator $L_1$  to form a Riesz basis.
\item[(c)]
 If the conditions of the Corollary 2 hold, then a normal EAF system
 of the operator $L_{1}$ forms a Riesz basis.
\end{itemize}
\end{theorem}

 The authors thank  Prof. B.S.Mitjagin and R.O.Hryniv
 for useful remarks.

\bigskip

A.A Shkalikov: Department of Mechanics and Mathematics, Lomonosov
Moscow State   University,

\noindent Leninskie Gory, Moscow, Russia

\medskip

\noindent e-mail: ashkalikov@yahoo.com

 \medskip

O.A. Veliev: Department  of Mathematics, Dogus University,
Acibadem, 34722,

\noindent Kadikoy, Istanbul, Turkey

\noindent e-mail:
oveliev@dogus.edu.tr

\end{document}